\documentclass[10pt]{article}
\usepackage{mathrsfs}
\usepackage{amsthm}
\usepackage{amsmath}
\usepackage{graphicx}
\usepackage{color}
\usepackage{amsfonts}
\newtheorem{theorem}{Theorem}[section]

\newtheorem{lemma}[theorem]{Lemma}

\newtheorem{corollary}[theorem]{Corollary}
\numberwithin{equation}{section}
\normalsize

\begin{document}
\title{\textbf{Asymptotic Behavior of Critical Infection Rates for Threshold-one Contact Processes on Lattices and Regular Trees}}
\author{Xiaofeng Xue \thanks{\textbf{E-mail}: masonxuexf@math.pku.edu.cn \textbf{Address}: School of Mathematical Sciences, Peking University, Beijing 100871, China.}\\ Peking University}

\date{}
\maketitle

\noindent {\bf Abstract}
In this paper we study threshold-one contact processes on lattices and regular trees. The asymptotic behavior of the critical infection rates as the degrees of the graphs growing to infinity are obtained. Defining \(\lambda_c\) as the supremum of infection rates which causes extinction of the process at equilibrium, we prove that \(n\lambda_c^{\textbf{T}^n}\rightarrow1\) and \(2d\lambda_c^{\textbf{Z}^d}\rightarrow1\) as \(n,d\rightarrow+\infty\). Our result is a development of the conclusion that \(\lambda_c^{\textbf{Z}^d}\leq\frac{2.18}{d}\) shown in \cite{Dur1991}. To prove our main result, a crucial lemma about the probability of a simple random walk on a lattice returning to zero is obtained. In details, the lemma is that \(\lim_{d\rightarrow+\infty}2dP\big(\exists n\geq1, S_n^{(d)}=0\big)=1\), where \(S_n^{(d)}\) is a simple random walk on \(\textbf{Z}^d\) with \(S_0^{(d)}=0\).


\section{Introduction}

In this paper we study threshold-one contact processes on lattices and regular trees. For a graph \(G\), threshold-one contact process on \(G\)
is with state space \(\{0,1\}^G\), which means that at each vertex on \(G\), there is a spin taking value \(0\) or \(1\). For each \(x\in G\) and any configuration \(\eta\in \{0,1\}^G\), we denote by \(\eta(x)\) the value of \(x\). For any \(t>0\), the configuration of the process at \(t\) is denoted by \(\eta_t\).
For each \(x\in G\) and \(t>0\), we define
\[
\eta_{t-}(x)=\lim_{s<t,s\uparrow t}\eta_s(x)
\]
as the value of \(x\) at the moment just before \(t\). For any \(x,y\in G\), we say that they are neighbors if there is an edge connecting them, denoted by \(x\sim y\).

Now we explain how the process evolves.  At the beginning, each vertex takes \(0\) or \(1\) according to some probability distribution. Then, the process evolves depending on independent Poisson processes \(\{N_x(t)\}_{x\in G}\) and  \(\{Y_x(t)\}_{x\in G}\). For each vertex \(x\in G\), \(N_x(\cdot)\) is with rate one while \(Y_x(\cdot)\) is with rate \(\lambda>0\). \(\lambda\) is called the infection rate. The value of \(x\) may flip only at the event times of \(N_x(\cdot)\) and \(Y_x(\cdot)\). For any event time \(s\) of \(N_x(\cdot)\), \(\eta_s(x)\) takes \(0\) no matter whatever \(\eta_{s-}(x)\) is. For any event time \(r\) of \(Y_x(\cdot)\), \(\eta_r(x)\) does not flip when \(\eta_{r-}(x)=1\). When \(\eta_{r-}(x)=0\), \(\eta_r(x)\) flips to \(1\) at \(r\) if and only if there exists a neighbor \(y\) of \(x\) such that \(\eta_{r-}(y)=1\). Therefore the threshold-one contact process \(\eta_t\) is a spin system (See the definition of spin systems in Chapter 3 of \cite{LIG1985}.) with flip rates given by
\begin{equation} \label{flip rate of Threshold-one contact process}
c(x,\eta)=
\begin{cases}
1 & \text{\quad if \quad} \eta(x)=1,\\
\lambda & \text{\quad if \quad} \eta(x)=0 \text{\quad and \quad} \sum_{y:y\sim x}\eta(y)>0,\\
0 & \text{otherwise}
\end{cases}
\end{equation}
for any configuration \(\eta\in \{0,1\}^G\).

Intuitively, the process describes the spread of an infection disease on a network. Each \(x\in G\) stands for an individual who may be infected by the disease. \(1\) and \(0\) represent the state `infected' and `healthy' respectively. An individual in the infected state will wait for an exponential time with rate one to be healed. An healthy individual will wait for an exponential time with rate \(\lambda\) to be infected if and only if there is at least one neighbor of it is in the infected state.

In later sections, we write \(\eta_t\) as \(\eta_t^\eta\) when \(\eta_0=\eta\in \{0,1\}^G\). We denote by \(\delta_1\) the configuration that all vertices take value \(1\). Since the threshold-one contact process is attractive (See the definition of `attractive' in Chapter 3 of \cite{LIG1985}.), it is easy to see that
\[
P\big(\eta_t^{\delta_1}(x)=1\big)
\]
decreases with \(t\) for each \(x\in G\). Hence it is reasonable to define
\[
\mu(x)=\lim_{t\rightarrow+\infty}P\big(\eta_t^{\delta_1}(x)=1\big)
\]
for each \(x\in G\). To distinguish processes on different graphs with different infection rates, we write \(\mu\) as \(\mu^G_{\lambda}\). According to the basic coupling of spin systems (See Chapter 3 of \cite{LIG1985}.), it is easy to see that
\[
\mu^G_{\lambda_1}(x)\geq \mu^G_{\lambda_2}(x)
\]
for \(\lambda_1>\lambda_2\). Therefore it is reasonable to define
\begin{equation}\label{defintion of critical value of threshold-one contact}
\lambda_c^G=\sup\{\lambda: \sup_{x\in G}\mu^G_{\lambda}(x)=0\}.
\end{equation}
\(\lambda_c^G\) is called the critical value of the infection rate. According to \eqref{defintion of critical value of threshold-one contact}, when \(\lambda<\lambda_c\), \(\eta_t\) converges weakly to \(\delta_0\), the configuration that all vertices take \(0\). Hence the disease is extinct when
\(\lambda<\lambda_c\). In this paper, we are concerned with the estimation of \(\lambda_c^G\) for \(G\) is a lattice or a regular tree. Our results will be introduced in following sections.

The threshold contact process is introduced in \cite{Dur1991} as a tool to study threshold voter model since when infection rate \(\lambda=1\), threshold voter models can be bounded below by threshold contact processes (See \cite{And1992}, \cite{Dur1991}, \cite{Han1999}, \cite{LIG1994}, \cite{LIG1999}, \cite{Xue2012}.). In \cite{Dur1991}, the threshold is considered to be one. It is shown in \cite{Dur1991} that threshold-one contact process has an additive dual process. Due to the additivity of the dual process, it is suggested that the threshold-one contact process has similar features with that of linear contact process which is additive and self-dual (See Chapter 6 of \cite{LIG1985}.). In recent years, more works are concerned on the case that the threshold is bigger than one such as \cite{Fonte2008}, \cite{Mou2009} and \cite{Xue2012}. It is studied in \cite{Fonte2008} and \cite{Xue2012} the critical infection rates and critical density points for threshold contact processes and threshold voter models on regular trees. It is showed in \cite{Dur1991} and \cite{Mou2009} that the critical infection rate for threshold contact process on lattice converges to \(0\) as the degree grows to infinity. This paper is a development of this result in the case of threshold one, as we give the asymptotic behavior of the critical infection rate.
\section{Main results}
\quad Now we introduce our main results. We obtain the asymptotic behavior of \(\lambda_c\) for the process on lattice and regular tree as the degree of the graph grows to infinity. In this paper, high-degree lattice with degree \(2d\) is denoted by \(\textbf{Z}^d\) while regular tree with degree \(n+1\) is denoted by \(\textbf{T}^n\). The following theorem is our main result.
\begin{theorem}\label{Theorem of asymptotic behavior of critical value}
\(\lambda_c^G\) is defined as that in \eqref{defintion of critical value of threshold-one contact}, then
\begin{equation}\label{critical value of lattice}
\lim_{d\rightarrow+\infty}2d\lambda_c^{\textbf{Z}^d}=1
\end{equation}
and
\begin{equation}\label{critical value of regular tree}
\lim_{n\rightarrow+\infty}n\lambda_c^{\textbf{T}^n}=1.
\end{equation}
\end{theorem}
Theorem \ref{Theorem of asymptotic behavior of critical value} shows that for lattices and regular trees, when the degree is large, \(\lambda_c\) is approximate to the reciprocal for the degree. In \cite{Dur1991}, Cox and Durrett shows that \(\lambda_c^{\textbf{Z}^d}\leq 2.18/d\). \eqref{critical value of lattice} is a development of this result. For classical linear contact process, similar asymptotic behaviors of critical value as  \eqref{critical value of lattice} and \eqref{critical value of regular tree} were proven in \cite{Grif1983} and \cite{Pem1992}. Theorem \ref{Theorem of asymptotic behavior of critical value} shows that the critical value of threshold-one contact process is with the same asymptotic behavior as that of linear contact process.

It is shown in \cite{Mou2009} that \(\lim_{d\rightarrow+\infty}\lambda_c^{\textbf{Z}^d}(K)=0\) where \(\lambda_c^{\textbf{Z}^d}(K)\) is the critical
value of threshold \(K\geq2\) contact process on \(\textbf{Z}^d\). By Theorem \ref{Theorem of asymptotic behavior of critical value}, it is natural to guess that
\begin{equation}
\lim_{d\rightarrow+\infty}2d\lambda_c^{\textbf{Z}^d}(K)=K
\end{equation}
for \(K\geq2\). But we have no idea whether this conjecture is right.

We divide the proof of Theorem \ref{Theorem of asymptotic behavior of critical value} into several sections. In Section \ref{section of lower bound}, we will prove \(\liminf_{d\rightarrow+\infty}2d\lambda_c^{\textbf{Z}^d}\geq1\) and \(\liminf_{n\rightarrow+\infty}n\lambda_c^{\textbf{T}^n}\geq1\) by giving a lower bound of \(\lambda_c\). In Section \ref{section of upper bound for tree}, we will give an upper bound of \(\lambda_c^{\textbf{T}^n}\) to accomplish the proof of \eqref{critical value of regular tree}. In Section \ref{section of upper bound for lattice}, we will give an upper bound of \(\lambda_c^{\textbf{Z}^d}\) to accomplish the proof of \eqref{critical value of lattice}.
\section{Lower bound}\label{section of lower bound}
\quad In this section we will give a lower bound of \(\lambda_c\). To do so, we introduce another stochastic process as a tool, which is denoted by \(\xi_t\). The state space of \(\xi_t\) on graph \(G\) is \(\textbf{N}^{G}\), where \(\textbf{N}\) is the set of nonnegative integers, which means that at each vertex there is an nonnegative integer. \(\xi_t\) evolves as following. \(\{N_x(\cdot)\}_{x\in G}\) and \(\{Y_x(\cdot)\}_{x\in G}\) are Poisson processes as that in the definition of threshold-one contact process. For each \(x\in G\), \(\xi(x)\) may change only at event times of \(N_x(\cdot)\) and \(Y_x(\cdot)\). At any event time \(s\) of \(N_x(\cdot)\), \(\xi_s(x)\) takes \(0\) no matter whatever \(\xi_{s-}(x)\) is. At any event time \(r\) of \(Y_x(\cdot)\), \(\xi(x)\) flips to \(\xi_r(x)=\xi_{r-}(x)+\sum_{y:y\sim x}\xi_{r-}(y)\) from \(\xi_{r-}(x)\). From the definition, it is easy to see that \(\xi_t\) is a linear model (See Chapter 9 of \cite{LIG1985}.). As a Markov process, \(\xi_t\) can also be described via its generator. For any \(\xi\in \textbf{N}^G\), \(x\in G\) and \(m\in \textbf{N}\), we define \(\xi^{x,m}\in \textbf{N}^G\) as
\begin{equation}
\xi^{x,m}(y)=
\begin{cases}
\xi(y) & \text{if\quad}y\neq x,\\
m & \text{if\quad}y=x.
\end{cases}
\end{equation}
Then the generator \(\Omega\) of \(\xi_t\) is given by
\begin{align}\label{generator of xi}
\Omega f(\xi)=&\sum_{x\in G}\big[f(\xi^{x,0})-f(\xi)\big]\\
&+\sum_{x\in G}\lambda\big[f(\xi^{x,\xi(x)+\sum_{y:y\sim x}\xi(y)})-f(\xi)\big]\notag
\end{align}
for any \(f\in C(\textbf{N}^{G})\) properly fast decaying.

Intuitively the process \(\xi_t\) counts the (degree of) seriousness of the disease throughout the process. At event times of \(Y_x(\cdot)\) an infected individual \(x\) is able to further infected by its neighbors. Whenever that occurs, we simply add the seriousness of the disease of \(x\) by the sum of all the seriousness of the disease of \(x\)'s neighbors.

We explain the connection between \(\xi_t\) and the threshold-one contact process \(\eta_t\). For each \(x\in G\) and \(t\geq0\), let \(\widehat{\eta}_t(x)=1_{\{\xi_t(x)>0\}}\). We claim that \(\widehat{\eta}_t\) is threshold-one contact process with flip rates given by \eqref{flip rate of Threshold-one contact process}. According to the definition of \(\xi_t\), at any event time \(s\) of \(N_x(\cdot)\), \(\xi_s(x)\) takes \(0\) and hence \(\widehat{\eta}_s(x)=0\). At any event time \(r\) of \(Y_x(\cdot)\), if \(\xi_{r-}(x)>0\), then \(\xi_r(x)\geq\xi_{r-}(x)>0\) and hence \(\widehat{\eta}(x)\) does not flip from \(1\) at \(r\). If \(\xi_{r-}(x)=\widehat{\eta}_{r-}(x)=0\), then \(\xi_{r}(x)>0\) if and only if there exists a neighbor \(y\) of \(x\) such that \(\xi_{r-}(y)>0\). In other words, \(\widehat{\eta}(x)\) flips from \(0\) to \(1\) at \(r\) if and only if there exists a neighbor \(y\) of \(x\) such that \(\widehat{\eta}_{r-}(y)=1\). Therefore \(\widehat{\eta}_t\) evolves as a threshold-one contact process. As a result, \(\xi_t^{\delta_1}\) and \(\eta_t^{\delta_1}\) with same infection rate \(\lambda\) on \(G\) can be coupled such that
\begin{equation}\label{xi_t control eta_t}
\eta_t(x)=1_{\{\xi_t(x)>0\}}
\end{equation}
for each \(x\in G\). By \eqref{xi_t control eta_t} and Chebyshev's inequality,
\begin{equation}\label{Exi bigger than Peta}
P\big(\eta_t^{\delta_1}(x)=1\big)=P\big(\xi_t^{\delta_1}(x)\geq1\big)\leq E\xi_t^{\delta_1}(x).
\end{equation}

Now we give a lower bound of \(\lambda_c^G\) where \(G\) is a simple regular graph. A regular graph is a graph where each vertex has
the same degree, and simple graphs are those containing no self-loops or multiple edges. \(\textbf{Z}^d\) and \(\textbf{T}^n\) are all simple regular graphs. The following theorem gives a lower bound of \(\lambda_c\).
\begin{theorem}\label{lower bound for regular graph}
For a simple regular graph \(G\) with degree \(r\),
\begin{equation}
\lambda_c^G\geq\frac{1}{r}.
\end{equation}
\end{theorem}
The following proposition is a direct corollary of Theorem \ref{lower bound for regular graph}.
\begin{corollary}\label{lower bound for lattice and regular tree}
\[
\lambda_c^{\textbf{Z}^d}\geq\frac{1}{2d}
\text{\quad and \quad}
\lambda_c^{\textbf{T}^n}\geq\frac{1}{n+1}.
\]
Hence,
\[
\liminf_{d\rightarrow+\infty}2d\lambda_c^{\textbf{Z}^d}\geq1
\text{\quad
and
\quad}
\liminf_{n\rightarrow+\infty}n\lambda_c^{\textbf{T}^n}\geq1
.\]
\end{corollary}
\proof[Proof of Theorem \ref{lower bound for regular graph}]
Consider \(\xi_t\) with infection rate \(\lambda\) on \(G\).
According to the generator \(\Omega\) of \(\xi_t\) given in \eqref{generator of xi}, we can prove that
\begin{equation}\label{ODE of xi}
\frac{d}{dt}E\xi_t^{\delta_1}(x)=-E\xi_t^{\delta_1}(x)+\lambda\sum_{y:y\sim x}E\xi_t^{\delta_1}(y)
\end{equation}
for each \(x\in G\).

Intuitively, \eqref{ODE of xi} is with the form \(\frac{d}{dt}Ef(\xi_t)=E\Omega f(\xi_t)\) as an `application' of Hille-Yosida Theorem (See Theorem 1.2.9 of \cite{LIG1985}). However, the state space \(\textbf{N}^{G}\) of \(\xi_t\) is not compact, which does not satisfy the condition of Hille-Yosida Theorem. To prove \eqref{ODE of xi} rigorously, we need Theorem 1.27 in Chapter 9 of \cite{LIG1985}. For more details, see
Appendix \ref{Appendixtwo}.

It is easy to verify that
\[
F_t(x)=\exp\{t(\lambda r-1)\}
\]
for each \(x\) is a solution of ODE \eqref{ODE of xi} with \(F_0=\textbf{1}\). According to classical theory of functional analysis, it is easy to see that ODE \eqref{ODE of xi} with initial condition \(E\xi_0^{\delta_1}=\textbf{1}\) has an unique solution. Therefore,
\[
E\xi_t^{\delta_1}(x)=F_t(x)=\exp\{t(\lambda r-1)\}
\]
for \(t\geq0\).
By \eqref{Exi bigger than Peta}, when \(\lambda<\frac{1}{r}\),
\[
\mu^G_\lambda(x)=\lim_{t\rightarrow+\infty}P\big(\eta_t^{\delta_1}(x)=1\big)\leq\lim_{t\rightarrow+\infty}E\xi_t^{\delta_1}(x)
=\lim_{t\rightarrow+\infty}\exp\{t(\lambda r-1)\}=0
\]
for each \(x\in G\), and hence
\[
\lambda_c^G\geq\frac{1}{r}.
\]

\qed

Since \(\textbf{Z}^d\) is simple regular graph with degree \(2d\) and \(\textbf{T}^n\) is simple regular graph with degree \(n+1\), Corollary \ref{lower bound for lattice and regular tree} follows from Theorem \ref{lower bound for regular graph} directly.
\section{Upper bound: the case of regular trees}\label{section of upper bound for tree}
\quad In this section we will give an upper bound of \(\lambda_c^{\textbf{T}^n}\) and accomplish the proof of \eqref{critical value of regular tree}.
A dual process \(A_t\) introduced in \cite{Dur1991} is crucial for our approach. The process \(A_t\) on \(\textbf{T}^n\) is with state space
\[
2^{\textbf{T}^n}=\{B:B\subseteq \textbf{T}^n\}.
\]
For each \(x\in \textbf{T}^n\), \(N_x(\cdot)\) and \(Y_x(\cdot)\) are Poisson processes as that in the definition of threshold-one contact process \(\eta_t\). \(A_t\) evolves as following. For each \(x\in \textbf{T}^n\) and any event time \(s\) of \(N_x(\cdot)\), \(A_s=A_{s-} \setminus\{x\}\). At any event time \(r\) of \(Y_x(\cdot)\), \(A_r=A_{r-}\cup\{y:y\sim x\}\) if \(x\in A_{r-}\), otherwise \(A_r=A_{r-}\). We write \(A_t\) as \(A_t^A\) when \(A_0=A\subseteq \textbf{T}^n\), then it is shown in \cite{Dur1991} that
\begin{equation}\label{dual for eta and A}
P\big(\eta_t^{\delta_1}(x)=1\big)=P\big(A_t^{\{x\}}\neq\emptyset\big)
\end{equation}
for each \(x\in \textbf{T}^n\) (See a simple proof of \eqref{dual for eta and A} in Part Two of \cite{LIG1999}).

We introduce a branching process \(S_t\in 2^{\textbf{T}^n}\) to bound below the growth of \(A_t\). To introduce \(S_t\), \(\textbf{T}^n\) is considered as an oriented regular tree that for each \(x\in \textbf{T}^n\), one neighbor of \(x\) is its `farther' while the other \(n\) neighbors of \(x\) are its sons. We denote by \(x\rightarrow y\) that \(y\) is a son of \(x\). \(S_t\) are evolves as following. For each \(x\in \textbf{T}^n\) and any event time \(s\) of \(N_x(\cdot)\), \(S_s=S_{s-} \setminus\{x\}\). At any event time \(r\) of \(Y_x(\cdot)\), \(S_r=\big(S_{r-}\cup\{y:x\rightarrow y\}\big)\setminus\{x\}\) if \(x\in S_{r-}\), otherwise \(S_r=S_{r-}\). We write \(S_t\) as \(S_t^A\) when \(S_0=A\subseteq \textbf{T}^n\).

By basic coupling, it is easy to see that
\[
A_t^{\{x\}}\supseteq S_t^{\{x\}}
\]
for each \(x\in \textbf{T}^n\) and any \(t\geq0\). Therefore,
\begin{equation}\label{PAt bigger than PSt}
P\big(\eta_t^{\delta_1}(x)=1\big)=P\big(A_t^{\{x\}}\neq\emptyset\big)\geq P\big(S_t^{\{x\}}\neq\emptyset\big).
\end{equation}
According to the definition of \(S_t\), for each \(x\in S_t\), \(x\) will be replaced by \(n\) sons with probability \(\frac{\lambda}{\lambda+1}\) or be kicked out from \(S_t\) without `compensation' with probability \(\frac{1}{\lambda+1}\). Therefore \(S_t\) is a branching process with offspring distribution with mean
\[
\frac{n\lambda}{\lambda+1}.
\]
The following theorem gives an upper bound of \(\lambda_c^{\textbf{T}^n}\).
\begin{theorem}\label{theorem of upper bound for regular tree}
\begin{equation}
\lambda_c^{\textbf{T}^n}\leq\frac{1}{n-1}
\end{equation}
and hence
\begin{equation}
\limsup_{n\rightarrow+\infty}n\lambda_c^{\textbf{T}^n}\leq1.
\end{equation}
\end{theorem}
\proof
For \(\lambda>\frac{1}{n-1}\),
\[
\frac{\lambda}{\lambda+1}n>1
.\]
Therefore according to classical theorems of branching process (See Chapter 3 of \cite{Hof2012}.),
\[
P\big(S_t^{\{x\}}\neq\emptyset \text{\quad for any \(t\geq0\)}\big)>0
\]
when \(\lambda>\frac{1}{n-1}\). By \eqref{PAt bigger than PSt},
\begin{align}
\mu^{\textbf{T}^n}_\lambda(x)&=\lim_{t\rightarrow+\infty}P\big(\eta_t^{\delta_1}(x)=1\big)=\lim_{t\rightarrow+\infty}P\big(A_t^{\{x\}}\neq\emptyset\big)\\
&\geq \lim_{t\rightarrow+\infty}P\big(S_t^{\{x\}}\neq\emptyset\big)=
P\big(S_t^{\{x\}}\neq\emptyset \text{\quad for any \(t\geq0\)}\big).\notag
\end{align}
Therefore, \(\mu^{\textbf{T}^n}_\lambda(x)>0\) for \(\lambda>\frac{1}{n-1}\) and
\[
\lambda_c^{\textbf{T}^n}\leq\frac{1}{n-1}
.\]

\qed

\eqref{critical value of regular tree} is a direct corollary of Corollary \ref{lower bound for lattice and regular tree} and Theorem \ref{theorem of upper bound for regular tree}. Furthermore, these two theorems show that
\[
\frac{1}{n+1}\leq\lambda_c^{\textbf{T}^n}\leq\frac{1}{n-1}
\]
for each \(n\geq1\).

\section{Upper bound: the case of lattices}\label{section of upper bound for lattice}
\quad In this section we will give an upper bound of \(\lambda_c^{\textbf{Z}^d}\) and accomplish the proof of \eqref{critical value of lattice}. The approach in Section \ref{section of upper bound for tree} fails here because there are many graph-loops on \(\textbf{Z}^d\) so that \(A_t\) can not be bounded below by a branching process. We are inspired a lot by the approach in Chapter 9 of \cite{LIG1985}. For some linear systems, the approach shows that the process is survival when the second moments are uniformly bounded.

As a tool, we introduce a stochastic process \(\zeta_t\) which is a modification of \(\xi_t\) introduced in Section \ref{section of lower bound}. For \(\zeta_t\) on \(\textbf{Z}^d\), the state space of \(\zeta_t\) is \([0,+\infty)^{\textbf{Z}^d}\), which means that at each vertex \(x\in \textbf{Z}^d\) there is a nonnegative real number. At event times of \(N_x(\cdot)\) and \(Y_x(\cdot)\), \(\zeta_t\) flips in the same way as that of \(\xi_t\). What different from \(\xi_t\) is that \(\zeta_t\) evolves according to an linear ODE between event times of Poisson processes.  In detail, for each \(x\in \textbf{Z}^d\), \(\zeta_s(x)=0\) at event time \(s\) of \(N_x(\cdot)\) while \(\zeta_r(x)=\zeta_{r-}(x)+\sum_{y:y\sim x}\zeta_{r-}(y)\) at event time \(r\) of \(Y_x(\cdot)\). Between any two adjacent event times of the Poisson processes \(N_x(\cdot)\) and \(Y_x(\cdot)\), \(\zeta_t(x)\) evolves according to the following ODE
\[
\frac{d}{dt}\zeta_t(x)=(1-2\lambda d)\zeta_t(x).
\]
The generator of \(\zeta_t\) is given by
\begin{align}\label{generator of zeta}
\Omega f(\zeta)=&\sum_{x\in \textbf{Z}^d}\big[f(\zeta^{x,0})-f(\zeta)\big]\\
&+\sum_{x\in \textbf{Z}^d}\lambda\big[f(\zeta^{x,\zeta(x)+\sum_{y:y\sim x}\zeta(y)})-f(\zeta)\big]+\sum_{x\in \textbf{Z}^d}f^\prime_x(\zeta)(1-2\lambda d)\zeta(x)\notag
\end{align}
for any \(\zeta\in [0,+\infty)^{\textbf{Z}^d}\) and \(f\in C^1\big([0,+\infty)^{\textbf{Z}^d}\big)\), where \(f^\prime_x(\zeta)\) is the partial derivative of \(f(\zeta)\) with respect to the coordinate \(\zeta(x)\) (See Chapter 9 of \cite{LIG1985} for more about generator of a linear system.).

The following Lemma shows that uniformly bounded second moments of \(\zeta_t\) ensure the survival of \(\eta_t\), which is crucial for our approach.
\begin{lemma}\label{second moment control survive}
If \(\lambda\) makes
\[
\sup_{t\geq0}E(\zeta_t^{\delta_1}(x))^2<+\infty
\]
for each \(x\in \textbf{Z}^d\),
then \(\lambda_c^{\textbf{Z}^d}\leq\lambda\).
\end{lemma}
Notice that \(E(\zeta_t^{\delta_1}(x))^2\) does not depending on \(x\) since \(\textbf{Z}^d\) is symmetric.
\proof[Proof of Lemma \ref{second moment control survive}]
For each \(x \in \textbf{Z}^d\) and \(t\geq0\), let \(\widetilde{\eta}_t(x)=1_{\{\zeta_t(x)>0\}}\). After a similar discussion with that of \(\xi_t\), it is easy to see that \(\widetilde{\eta}_t\) is also a threshold-one contact process with flip rates given by \eqref{flip rate of Threshold-one contact process}. Therefore, \(\eta_t^{\delta_1}\) and \(\zeta_t^{\delta_1}\) with same infection rate \(\lambda\) on \(\textbf{Z}^d\) can be coupled such that
\[
\eta_t^{\delta_1}(x)=1_{\{\zeta_t^{\delta_1}(x)>0\}}
\]
for each \(x\in \textbf{Z}^d\). Then by H\"{o}lder inequality,
\begin{equation}
P\big(\eta_t^{\delta_1}(x)=1\big)=P\big(\zeta_t^{\delta_1}(x)>0\big)\geq\frac{\big(E\zeta^{\delta_1}_t(x)\big)^2}{E\big(\zeta_t^{\delta_1}(x)\big)^2}
.\end{equation}

According to the generator of \(\zeta_t\) given in \eqref{generator of zeta} and a similar proof with that of \eqref{ODE of xi} (See Appendix \ref{Appendixtwo}.),
\begin{align}\label{ODE of zeta}
\frac{d}{dt}E\zeta^{\delta_1}_t(x)=-E\zeta_t^{\delta_1}(x)+\lambda\sum_{y:y\sim x}E\zeta_t^{\delta_1}(y)+(1-2\lambda d)E\zeta_t^{\delta_1}(x)
\end{align}
for each \(x\in \textbf{Z}^d\). Due to the symmetry of \(\textbf{Z}^d\), \(E\zeta^{\delta_1}_t(x)\) does not depending on \(x\). Therefore,
\[
\frac{d}{dt}E\zeta^{\delta_1}_t(x)=-E\zeta_t^{\delta_1}(x)+2d\lambda E\zeta_t^{\delta_1}(x)+(1-2\lambda d)E\zeta_t^{\delta_1}(x)=0
\]
and
\[
E\zeta^{\delta_1}_t(x)\equiv E\zeta^{\delta_1}_0(x)=1
\]
for any \(t\geq0\).
Therefore, \(P\big(\eta_t^{\delta_1}(x)=1\big)\geq 1/E\big(\zeta_t^{\delta_1}(x)\big)^2\). If \(\lambda\) makes
\[
\sup_{t\geq0}E(\zeta_t^{\delta_1}(x))^2<+\infty
\]
then
\[
\mu^{\textbf{Z}^d}_\lambda(x)=\lim_{t\rightarrow+\infty}P\big(\eta_t^{\delta_1}(x)=1\big)\geq 1/\sup_{t\geq0}E(\zeta_t^{\delta_1}(x))^2>0
\]
and hence \(\lambda\geq\lambda_c^{\textbf{Z}^d}\).

\qed

Now the main problem is to find \(\lambda\) making \(\sup_{t\geq0}E(\zeta_t^{\delta_1}(x))^2<+\infty\). First we give the ODE which \(\{E(\zeta_t^{\delta_1}(x))^2\}_{x\in \textbf{Z}^d}\) satisfying. By the symmetry of \(\textbf{Z}^d\), we define
\[
G_t(x)=E\big[\zeta_t^{\delta_1}(0)\zeta_t^{\delta_1}(x)\big]=E\big[\zeta_t^{\delta_1}(y)\zeta_t^{\delta_1}(x+y)\big]
\]
for any \(x,y\in \textbf{Z}^d\) and \(t\geq0\). Then \(E(\zeta_t^{\delta_1}(x))^2=G_t(0)\).
According to the generator of \(\zeta_t\), we can show that \(\{G_t(x)\}_{x\in \textbf{Z}^d}\) satisfies
\begin{equation}\label{ODE of Gzero}
\frac{d}{dt}G_t(x)=4\lambda d\big[\frac{1}{2d}\sum_{y:y\sim x}G_t(y)-G_t(x)\big]
\end{equation}
for any \(x\neq0\) and
\begin{align}\label{ODE of Gnonzero}
\frac{d}{dt}G_t(0)&=(1-4\lambda d)G_t(0)+2\lambda\sum_{y:y\sim 0}G_t(y)+\lambda\sum_{y:y\sim 0}\sum_{z:z\sim 0}G_t(y+z)\\
&=(1-2\lambda d)G_t(0)+2\lambda\sum_{y:y\sim 0}G_t(y)+\lambda\sum_{y:y\sim 0}\sum_{z:z\sim 0,\atop z\neq-y}G_t(y+z).\notag
\end{align}
In other words,
\[
\frac{d}{dt}G_t=QG_t,
\]
where \(Q\) is a \(\textbf{Z}^d*\textbf{Z}^d\) matrix such that
\begin{equation}\label{Q}
Q(x_1,x_2)=
\begin{cases}
-4\lambda d & \text{\quad if \(x_1=x_2\neq0\)}\\
2\lambda    & \text{\quad if \(x_1\neq0, x_2\sim x_1\)}\\
1-2\lambda d & \text{\quad if \(x_1=x_2=0\)}\\
2\lambda  & \text{\quad if \(x_1=0, x_2\sim 0\)}\\
\sum_{(y,z):y\sim 0 \atop z\sim 0, y+z=x_2}\lambda & \text{\quad if \(x_1=0, \|x_2\|=2\)}\\
0 & \text{\quad else}
\end{cases}
\end{equation}
for \(x_1,x_2\in \textbf{Z}^d\).

\eqref{ODE of Gzero} and \eqref{ODE of Gnonzero} are also with the form \(\frac{d}{dt}Ef(\zeta_t)=E\Omega f(\zeta_t)\) as \eqref{ODE of xi}. To prove these two equations rigorously, we need Theorem 3.1 of Chapter 9 of \cite{LIG1985}. For more details, see Appendix \ref{Appendixtwo}.

The following Lemma gives a sufficient condition for \(\sup_{t\geq0}E(\zeta_t^{\delta_1}(x))^2<+\infty\).
\begin{lemma}\label{harmonic h control uniform bound second moment}
If there exists a function \(h:\textbf{Z}^d\rightarrow \textbf{R}\) such that
\begin{equation}\label{h bounded}
0<\inf_{x \in \textbf{Z}^d} h(x)\leq\sup_{x\in \textbf{Z}^d}h(x)<+\infty
\end{equation}
and
\begin{equation}\label{h Q-harmonic}
Qh=0,
\end{equation}
then
\[
\sup_{t\geq0}E(\zeta_t^{\delta_1}(x))^2\leq\frac{\sup_{x\in \textbf{Z}^d}h(x)}{\inf_{x \in \textbf{Z}^d} h(x)}<+\infty
.\]
\end{lemma}
The following proof of Lemma \ref{harmonic h control uniform bound second moment} need several characters of the matrix \(Q\). We will prove these characters rigorously in Appendix \ref{character of Q}.
\proof
In Theorem \ref{upper bound of Qelement} of Appendix \ref{character of Q} we will show that
\[
|Q^n(x,y)|\leq(1+8\lambda d+4\lambda d^2)^n
\]
for any \(x,y\in \textbf{Z}^d\).

Therefore it is reasonable to define
\[
\exp\{tQ\}=\sum_{n=0}^{+\infty}\frac{(tQ)^n}{n!}
\]
for any \(t\geq0\).

We denote by \(L^\infty(\textbf{Z}^d)\) the set of bounded functions on \(\textbf{Z}^d\) and define
\[
\|f\|_\infty=\sup_{x\in \textbf{Z}^d}|f(x)|
\]
for \(f\in L^\infty(\textbf{Z}^d)\).

According to classical theorems of linear ODE, we will show in Theorem \ref{unique solution of linear ode} of Appendix \ref{character of Q} that the unique solution to the following equation
\begin{align}
\frac{d}{dt}f_t=Qf_t
\end{align}
with initial condition \(f_0\in L^\infty(\textbf{Z}^d)\) is
\[
f_t=\exp\{tQ\}f_0
\] and \(\|f_t\|_\infty \leq \exp\{t(1+8\lambda d+4\lambda d^2)\}\|f_0\|_\infty\).

As a result,
\begin{align*}
E\big(\zeta_t^{\delta_1}(0)\big)^2&=G_t(0)=\sum_{x\in \textbf{Z}^d}\exp\{tQ\}(0,x)G_0(x)\\
&=\sum_{x\in \textbf{Z}^d}\exp\{tQ\}(0,x)E\big(\zeta_0^{\delta_1}(0)\zeta_0^{\delta_1}(x)\big)=\sum_{x\in \textbf{Z}^d}\exp\{tQ\}(0,x).
\end{align*}

According to the definition of \(\exp\{tQ\}\) and Fubini Theorem,
\[
Q\exp\{tQ\}=\sum_{n=0}^{+\infty}\frac{t^nQ^{n+1}}{n!}=\exp\{tQ\}Q.
\]

Since \(Qh=0\),
\[
\frac{d}{dt}\exp\{tQ\}h=Q\exp\{tQ\}h=\exp\{tQ\}Qh=0
\]
and \(\exp\{tQ\}h\equiv h\) for \(t\geq0\).
Therefore,
\[
h(0)=\exp\{tQ\}h(0)=\sum_{x\in \textbf{Z}^d}\exp\{tQ\}(0,x)h(x)
\]
for \(t\geq0\).
In Theorem \ref{positive of Qelement} of Appendix \ref{character of Q} we will show that
\[
\exp\{tQ\}(x,y)\geq0
\]
for any \((x,y)\in \textbf{Z}^d\).
Then
\begin{align*}
E\big(\zeta_t^{\delta_1}(0)\big)^2&=\sum_{x\in \textbf{Z}^d}\exp\{tQ\}(0,x)\\
&\leq \sum_{x\in \textbf{Z}^d}\exp\{tQ\}(0,x)\frac{h(x)}{\inf_{y \in \textbf{Z}^d} h(y)}\\
&=\frac{h(0)}{\inf_{y \in \textbf{Z}^d} h(y)}\\
&\leq\frac{\sup_{x\in \textbf{Z}^d}h(x)}{\inf_{x \in \textbf{Z}^d} h(x)}
\end{align*}
for \(t\geq0\) and the proof complete.

\qed

To construct \(h\) satisfying \eqref{h Q-harmonic}, we consider simple random walk \(S_n^{(d)}\) on \(\textbf{Z}^d\).
Let \(\tau_0^{(d)}=\inf\{m\geq1: S_m^{(d)}=0\}\). For \(d\geq1\), we define \(F_d:\textbf{Z}^d\rightarrow \textbf{R}\) as
\[
F_d(x)=P(\tau_0^{(d)}<+\infty|S_0=x)
\]
for \(x\in Z^d\setminus\{0\}\) and \(F_d(0)=1\). Let \(e_1^{(d)}=(1,\underbrace{0,0,\ldots,0}_{d-1 \text{th}})\), then the following estimation
of \(F_d(e_1^{(d)})\) is crucial for us to construct \(h\).
\begin{lemma}
\begin{equation} \label{asymptotic behavior of hitting prob}
\lim_{d\rightarrow+\infty}2dF_d(e_1^{(d)})=1.
\end{equation}
\end{lemma}
We do not know whether \eqref{asymptotic behavior of hitting prob} has been proven in early references about simple random walk. We searched several famous books such as \cite{Law2010} and \cite{Spi1976} but can not find this conclusion, so we give our own proof of \eqref{asymptotic behavior of hitting prob} in Appendix \ref{Appendixone}.

By \eqref{asymptotic behavior of hitting prob},
\[
\frac{1}{4d\big[1-(d+1)F_d(e_1^{(d)})\big]}>0
\]
for sufficient large \(d\). Finally  we can construct \(h\) and give an upper bound of \(\lambda_c^{\textbf{Z}^d}\).
\begin{theorem}\label{theorem of upper bound for lattice}
For sufficient large \(d\) such that \(\frac{1}{4d\big[1-(d+1)F_d(e_1^{(d)})\big]}>0\) and
\[
\lambda>\frac{1}{4d\big[1-(d+1)F_d(e_1^{(d)})\big]}
,\]
we define
\[
b_\lambda=\frac{4d\lambda\big[1-(d+1)F_d(e_1^{(d)})\big]-1}{1+4d^2\lambda}
\]
and
\[
h(x)=F_d(x)+b_\lambda
\]
for each \(x\in \textbf{Z}^d\). Then \(h\) satisfies \eqref{h bounded} and \eqref{h Q-harmonic}. As a result,
\[
\lambda_c^{\textbf{Z}^d}\leq \frac{1}{4d\big[1-(d+1)F_d(e_1^{(d)})\big]}.
\]
and
\[
\limsup_{t\rightarrow+\infty}2d\lambda_c^{\textbf{Z}^d}\leq1
.\]
\end{theorem}
\proof
When \(\lambda>\frac{1}{4d\big[1-(d+1)F_d(e_1^{(d)})\big]}\),
\[
0<b_\lambda\leq \inf h(x)\leq\sup h(x)\leq1+b_\lambda<+\infty.
\]
Hence \(h\) satisfies \eqref{h bounded}. For \eqref{h Q-harmonic}, when \(x\neq0\),
\[
4\lambda d[\frac{1}{2d}\sum_{y:y\sim x}h(y)-h(x)]=4\lambda d[\frac{1}{2d}\sum_{y:y\sim x}F_d(y)-F_d(x)]=0
\]
according to the probability transition of \(S_n\). For the case of \(0\),
\begin{align*}
&(1-4\lambda d)h(0)+2\lambda\sum_{y:y\sim 0}h(y)+\lambda\sum_{y:y\sim0}\sum_{z:z\sim0}h(y+z)\\
=&(1+4\lambda d^2)b_\lambda+(1-4\lambda d)+2\lambda\sum_{y:y\sim0}F_d(y)+\lambda\sum_{y:y\sim0}\big[1+\sum_{z:z\sim0,\atop z\neq-y}F_d(y+z)\big]\\
=&(1+4\lambda d^2)b_\lambda+(1-4\lambda d)+2\lambda\sum_{y:y\sim0}F_d(y)+\lambda\sum_{y:y\sim0}2dF_d(y)\\
=&(1+4\lambda d^2)b_\lambda+(1-4\lambda d)+4\lambda d(d+1)F_d(e_1^{(d)})\\
=&(1+4\lambda d^2)b_\lambda+1-4\lambda d[1-(d+1)F_d(e_1^{(d)})]\\
=&0
\end{align*}
according to the definition of \(b_\lambda\) and \(h\).
Notice that during the calculation, we use that \(F_d(y)=F_d(e_1^{(d)})\) for \(y\sim 0\) since \(\textbf{Z}^d\) is symmetric and \[\frac{1}{2d}+\frac{1}{2d}\sum_{z:z\sim0,\atop z\neq-y}F_d(y+z)=F_d(y)\] for \(y\sim 0\).

The calculation above shows that \(Qh=0\). By Lemma \ref{second moment control survive} and Lemma \ref{harmonic h control uniform bound second moment}, \[
\lambda\geq\lambda_c^{\textbf{Z}^d}
\]
for any \(\lambda>\frac{1}{4d\big[1-(d+1)F_d(e_1^{(d)})\big]}\) and hence \(\lambda_c^{\textbf{Z}^d}\leq\frac{1}{4d\big[1-(d+1)F_d(e_1^{(d)})\big]}\). Furthermore, \(\limsup_{t\rightarrow+\infty}2d\lambda_c^{\textbf{Z}^d}\leq1\) holds by \eqref{asymptotic behavior of hitting prob}.

\qed

\eqref{critical value of lattice} is a direct corollary of Corollary \ref{lower bound for lattice and regular tree} and Theorem \ref{theorem of upper bound for lattice}. For large \(d\), we shows that
\begin{equation}
\frac{1}{2d}\leq\lambda_c^{\textbf{Z}^d}\leq\frac{1}{4d\big[1-(d+1)F_d(e_1^{(d)})\big]}
.
\end{equation}
Now the whole proof of Theorem \ref{Theorem of asymptotic behavior of critical value} is accomplished.

\appendix
\section{Appendix}
\subsection{Proof of \eqref{asymptotic behavior of hitting prob}}\label{Appendixone}
\proof[Proof of \eqref{asymptotic behavior of hitting prob}]
According to classical theory of simple random walk (See \cite{Law2010} and \cite{Spi1976}.),
\[
F_d(e_1^{(d)})=\frac{G_d(0,0)-1}{G_d(0,0)},
\]
where \(S_n^{(d)}\) is simple random walk on \(\textbf{Z}^d\) with \(S_0^{(d)}=0\) and
\[
G_d(0,0)=1+\sum_{n=1}^{+\infty}P(S^{(d)}_{2n}=0)
.\]
Hence we only need to show that \(\lim_{d\rightarrow+\infty}2d[G_d(0,0)-1]=1\).
\[
G_d(0,0)-1=2d(\frac{1}{2d})^2+\sum_{n=2}^{+\infty}P(S^{(d)}_{2n}=0)=\frac{1}{2d}+\sum_{n=2}^{+\infty}P(S^{(d)}_{2n}=0)
.\]
Hence we only need to show that \(\lim_{d\rightarrow+\infty}d\sum_{n=2}^{+\infty}P(S^{(d)}_{2n}=0)=0\).
Let
\[
H_d(1)=\sum_{n=2}^{d}P(S^{(d)}_{2n}=0)
\]
and
\[
H_d(k)=\sum_{n=(k-1)d+1}^{kd}P(S^{(d)}_{2n}=0)
\]
for \(k\geq2\). Then \(\sum_{n=2}^{+\infty}P(S^{(d)}_{2n}=0)=\sum_{k=1}^{+\infty} H_d(k)\).
\begin{align*}
H_d(1)&\leq\sum_{n=2}^d{2n \choose n}d^nn!(\frac{1}{2d})^{2n}\\
&=\sum_{n=2}^dL(n,d),
\end{align*}
where
\[
L(n,d)=\frac{(2n-1)!!}{(2d)^n}
.\]
\(L(n,d)=L(n-1,d)\frac{2n-1}{2d}\), hence \(L(n,d)\) decreases with \(n\) when \(n<\lceil d\rceil\) and increases with \(n\) when \(n\geq\lceil d\rceil\). When \(n\leq\frac{d}{2}\), \(L(n,d)\leq \frac{1}{2}L(n-1,d)\) and \(L(n,d)\leq \frac{1}{2^{n-2}}L(2,d)\).
Therefore,
\begin{align}
H_d(1)&\leq L(2,d)\sum_{n=2}^{\lfloor\frac{d}{2}\rfloor}\frac{1}{2^{n-2}}+(d-\lfloor\frac{d}{2}\rfloor)L(\lfloor\frac{d}{2}\rfloor,d)\\
&\leq \frac{3}{2d^2}+dL(\lfloor\frac{d}{2}\rfloor,d)\notag.
\end{align}
By Stirling formula,
\begin{align}
L(\lfloor\frac{d}{2}\rfloor,d)&=\frac{(2\lfloor \frac{d}{2}\rfloor)!}{(2d)^{\lfloor \frac{d}{2}\rfloor}\lfloor \frac{d}{2}\rfloor!2^{\lfloor \frac{d}{2}\rfloor}}\notag\\
&=\frac{\sqrt{2(1+o(1))}}{(2e)^{\lfloor \frac{d}{2}\rfloor}}.
\end{align}
Therefore,
\[
H_d(1)\leq \frac{3}{2d^2}+\frac{2d}{(2e)^{\lfloor \frac{d}{2}\rfloor}}
\]
for sufficiently large \(d\) and \(\lim_{d\rightarrow+\infty}dH_d(1)=0\)

Since \(L(n,d)\) increases with \(n\) when \(n\geq d+1\),
\begin{align*}
H_d(2)&\leq\sum_{n=d+1}^{2d}L(n,d)\\
&\leq dL(2d,d)\\
&=d\sqrt{2(1+o(1))}\big(\frac{2}{e}\big)^{2d}
\end{align*}
by Stirling formula.
Therefore \(H_d(2)\leq2d\big(\frac{2}{e}\big)^{2d}\) for sufficiently large \(d\) and \(\lim_{d\rightarrow+\infty}dH_d(2)=0\).

Finally we will show that \(\lim_{d\rightarrow+\infty}d\sum_{k=3}^{+\infty}H_d(k)=0\). For \(n\geq1\), we define
\[
\beta(n)=\frac{n!}{\sqrt{2\pi n}\big(\frac{n}{e}\big)^n}.
\]
By Stirling formula, \(\lim_{n\rightarrow +\infty}\beta(n)=1\). Hence there exists \(N_1\) such that \(\frac{\beta_{2n}}{\beta_{n}}<2\) for any \(n\geq N_1\). For \(1\leq j\leq d\) and \(k\geq3\),
\begin{align*}
P(S_{2(kd+j)}^{(d)}=0)&=\sum_{l_1+l_2+\ldots+l_d=kd+j}\frac{\big[2(kd+j)\big]!}{(l_1!)^2(l_2!)^2\ldots(l_d!)^2}\big(\frac{1}{2d}\big)^{2(kd+j)}\\
&=\frac{\big[2(kd+j)\big]!}{(kd+j)!(kd+j)!}\sum_{l_1+l_2+\ldots+l_d=kd+j}\Big[\frac{(kd+j)!}{l_1!l_2!\ldots l_d!}\Big]^2\big(\frac{1}{2d}\big)^{2(kd+j)}
\end{align*}
Since \(l_1!l_2!\ldots l_d!\geq(k!)^{d-j}((k+1)!)^j\) and
\[
\sum_{l_1+l_2+\ldots+l_d=kd+j}\frac{(kd+j)!}{l_1!l_2!\ldots l_d!}=d^{kd+j},
\]
we have
\begin{align*}
&P(S_{2(kd+j)}^{(d)}=0)\\
&\leq\frac{\beta(2(kd+j))2^{2(kd+j)}}{\beta^2(kd+j)}\sqrt{\frac{1}{\pi(kd+j)}}
\frac{(kd+j)!d^{kd+j}}{(k!)^{d-j}((k+1)!)^j}\big(\frac{1}{2d}\big)^{2(kd+j)}\\
&=\sqrt{2}\frac{\beta(2(kd+j))}{\beta(kd+j)}\frac{(k+\frac{j}{d})^{kd+j}}{e^{kd+j}(k!)^{d-j}((k+1)!)^j}\\
&\leq\sqrt{2}\frac{\beta(2(kd+j))}{\beta(kd+j)}\Big[\frac{(k+1)^k}{e^kk!}\Big]^d.
\end{align*}
Let \(M_k=\frac{(k+1)^k}{e^kk!}\), then \(\frac{M_{k+1}}{M_k}=(1+\frac{1}{k+1})^{k+1}/e<1\), and hence \(\sup_{k\geq2}M_k=M_2=\frac{9}{2e^2}<1\).
By Stirling Formula, \(\lim_{k\rightarrow+\infty}\sqrt{2\pi k}M_k=e\) and hence \(C=\sup_{k\geq2}\sqrt{k}M_k<+\infty\). Choose \(N_2\) such that
\(\frac{C}{\sqrt{N_2}}<\frac{1}{2}\), then for \(d>N_1\),
\begin{align*}
\sum_{k=3}^{+\infty}H_d(k)&\leq\sum_{k=3}^{+\infty}2\sqrt{2}dM_{k-1}^d\\
&\leq 2\sqrt{2}dN_2M_2^d+2\sqrt{2}dC^d\sum_{m=N_2+1}^{+\infty}k^{-\frac{d}{2}}\\
&\leq 2\sqrt{2}dN_2M_2^d+2\sqrt{2}dC^d\int_{N_2}^{+\infty}x^{-\frac{d}{2}}dx\\
&=2\sqrt{2}dN_2M_2^d+\frac{4\sqrt{2}dN_2}{d-2}\big(\frac{C}{\sqrt{N_2}}\big)^d,
\end{align*}
and
\[
d\sum_{k=3}^{+\infty}H_d(k)\leq 2\sqrt{2}d^2N_2M_2^d+\frac{4\sqrt{2}d^2N_2}{d-2}(\frac{1}{2})^d.
\]
Since \(M_2<1\),
\[
\lim_{d\rightarrow+\infty}d\sum_{k=3}^{+\infty}H_d(k)=0
.\]
As a result,
\[
\lim_{d\rightarrow+\infty}d\sum_{n=2}^{+\infty}P(S^{(d)}_{2n}=0)=\lim_{d\rightarrow+\infty}d\big(H_d(1)+H_d(2)+\sum_{k=3}^{+\infty}H_d(k)\big)=0.
\]

\qed

\subsection{Proof of \eqref{ODE of xi}, \eqref{ODE of zeta}, \eqref{ODE of Gzero} and \eqref{ODE of Gnonzero}.}\label{Appendixtwo}
\quad In this subsection we give the rigorous proofs of \eqref{ODE of xi}, \eqref{ODE of zeta}, \eqref{ODE of Gzero} and \eqref{ODE of Gnonzero}.
\proof[Proof of \eqref{ODE of xi}]
For any \(t>0\), we define that \(\beta_t=\xi_{\frac{1}{1+\lambda}t}\). According to the flip-rates of \(\xi_t\), \(\beta_t\) is a standard linear system introduced in Chapter 9.0 of \cite{LIG1985} with
\[
a(u,v)=0
\]
for any \(u, v\in G\) and
\begin{equation}
A_x(u,v)=
\begin{cases}
1 & \text{\quad if \quad} u=v\neq x,\\
0 & \text{\quad else}
\end{cases}
\end{equation} with probability \(\frac{1}{1+\lambda}\) and
\begin{equation}
A_x(u,v)=
\begin{cases}
1 & \text{\quad if \quad } u=v,\\
1 & \text{\quad else if \quad} u=x \text{\quad and \quad} v\sim x,\\
0 & \text{\quad else}
\end{cases}
\end{equation} with probability \(\frac{\lambda}{1+\lambda}\)
for any \(u, v, x\in G\).

According to Theorem 1.27 of Chapter 9 of \cite{LIG1985},
\begin{equation}
\frac{d}{dt}E\beta_t^{\delta_1}(x)=\sum_y\gamma(x,y)E\beta_t^{\delta_1}(y)
\end{equation}
where
\begin{equation}
\gamma(x,y)=a(x,y)+
\begin{cases}
{\rm E}\sum_uA_u(x,y) & \text{\quad if \quad} x\neq y,\\
{\rm E}\sum_u[A_u(x,x)-1]  & \text{\quad else.}
\end{cases}
\end{equation}
By the definition of \(a(\cdot,\cdot)\), \(\{A_x(\cdot,\cdot)\}_{x\in G}\) and direct calculation,
\begin{equation}
\gamma(x,y)=
\begin{cases}
-\frac{1}{1+\lambda} & \text{\quad if \quad} x=y,\\
\frac{\lambda}{1+\lambda}& \text{\quad else if \quad} y\sim x,\\
0 & \text{\quad else}
\end{cases}
\end{equation}
and
\begin{equation}\label{ODE of beta}
\frac{d}{dt}E\beta_t^{\delta_1}(x)=-\frac{1}{1+\lambda}E\beta_t^{\delta_1}(x)+\sum_{y:y\sim x}\frac{\lambda}{1+\lambda}E\beta_t^{\delta_1}(y).
\end{equation}
Since \(E\xi_t^{\delta_1}(x)=E\beta_{(1+\lambda)t}^{\delta_1}(x)\) for any \(x\in G\), \eqref{ODE of xi} is a direct corollary of \eqref{ODE of beta}.

\qed

\proof[Proof of \eqref{ODE of zeta}]
For any \(t>0\), we define that \(\alpha_t=\zeta_{\frac{1}{1+\lambda}t}\). Then \(\alpha_t\) is a standard linear model with
\[
a(x,y)=
\begin{cases}
\frac{1-2d\lambda}{1+\lambda} & \text{\quad if \quad} x=y,\\
0 & \text{\quad else} \\
\end{cases}
\]
for any \(x,y\in \textbf{Z}^d\) and the same \(\{A_x(\cdot,\cdot)\}_{x\in \textbf{Z}^d}\) as that of \(\beta_t\) in the proof of \eqref{ODE of xi}.
As we have done in the proof of \eqref{ODE of xi}, we can obtain \eqref{ODE of zeta} by directly applying Theorem 1.27 of Chapter 9 of \cite{LIG1985}. We omit the details.

\qed

\proof[Proof of \eqref{ODE of Gzero} and \eqref{ODE of Gnonzero}]
\(\alpha_t\) is the same as that in the proof of \eqref{ODE of zeta}. We use \(g(t,x,y)\) to denote \(E\alpha_t^{\delta_1}(x)\alpha_t^{\delta_1}(y)\) for any \(t\geq0\) and \(x,y\in \textbf{Z}^d\). \(\{A_x(\cdot,\cdot)\}_{x\in \textbf{Z}^d}\) and \(a(\cdot,\cdot)\) are the same as that in the proof of \eqref{ODE of zeta}. By direct calculation it is easy to verify that
\[
{\rm E}\sum_{x\in \textbf{Z}^d}\big[|A_x(u,u)-1|+\sum_{v:v\neq u}A_x(u,v)\big]^2<+\infty
\] for any \(u\in \textbf{Z}^d\).

Then according to Theorem 3.1 of Chapter 9 of \cite{LIG1985}, \(g(t,x,y)\) satisfies that \(g(0,x,y)=1\) and
\begin{equation}\label{ODE of alpha}
\frac{d}{dt}g(t,x,y)=\sum_{u,v\in \textbf{Z}^d}q\big((x,y),(u,v)\big)g(t,u,v)
\end{equation}
where
\[
q\big((x,y),(u,v)\big)=
\begin{cases}
{\rm E}\sum_zA_z(x,u)A_z(y,v) & \text{\quad if \quad} u\neq x, v\neq y,\\
{\rm E}\sum_zA_z(x,x)A_z(y,v)+a(y,v) & \text{\quad else if \quad} u=x, v\neq y,\\
{\rm E}\sum_zA_z(x,u)A_z(y,y)+a(x,u) & \text{\quad else if \quad} u\neq x, v\neq y,\\
{\rm E}\sum_z\big[A_z(x,x)A_z(y,y)-1\big]+a(x,x)+a(y,y) & \text{\quad else.}
\end{cases}
\]
According to the definition of \(\{A_x(\cdot,\cdot)\}_{x\in \textbf{Z}^d}\) and \(a(\cdot,\cdot)\) of \(\alpha_t\),
\begin{equation}\label{oneform of q}
q\big((x,x),(u,v)\big)=
\begin{cases}
\frac{1-4\lambda d}{1+\lambda} & \text{\quad if \quad} u=x, v=x,\\
\frac{\lambda}{1+\lambda} & \text{\quad else if \quad} u=x, v\sim x,\\
\frac{\lambda}{1+\lambda} & \text{\quad else if \quad} u\sim x, v=x,\\
\frac{\lambda}{1+\lambda} & \text{\quad else if \quad} u\sim x, v\sim x,\\
0 & \text{\quad else}\\
\end{cases}
\end{equation}
and
\begin{equation}\label{twoform of q}
q\big((x,y),(u,v)\big)=
\begin{cases}
-\frac{4\lambda d}{1+\lambda} & \text{\quad if \quad} u=x, v=y,\\
\frac{\lambda}{1+\lambda} & \text{\quad else if \quad} u=x, v\sim y,\\
\frac{\lambda}{1+\lambda} & \text{\quad else if \quad} u\sim x, v=y,\\
0 & \text{\quad else}
\end{cases}
\end{equation}
for \(x\neq y\).

Since \(G_t(x)=g\big((1+\lambda)t,0,x\big)=g\big((1+\lambda)t,y,y+x\big)\) for any \(x,y\in \textbf{Z}^d\), \eqref{ODE of Gzero} and \eqref{ODE of Gnonzero} are direct corollaries of \eqref{ODE of alpha}, \eqref{oneform of q} and \eqref{twoform of q}.

\qed

\subsection{Characters of \(Q\)}\label{character of Q}
\begin{theorem}\label{upper bound of Qelement}
\(Q\) is the same as that in \eqref{Q}. For any \(x,y\in \textbf{Z}^d\),
\begin{equation}\label{poly bound of Q}
|Q^n(x,y)|\leq(1+8\lambda d+4\lambda d^2)^n.
\end{equation}
\end{theorem}
\proof
By \eqref{Q} and direct calculation it is easy to see that \eqref{poly bound of Q} holds for \(n=1\).
When \eqref{poly bound of Q} holds for some \(n\geq1\),
\begin{align*}
|Q^{n+1}(x,y)|&=|QQ^n(x,y)|\\
&=|4\lambda d\big[\frac{1}{2d}\sum_{z:z\sim x}Q^n(z,y)-Q^n(x,y)\big]|\\
&\leq 8\lambda d(1+8\lambda d+4\lambda d^2)^{n}\\
&\leq (1+8\lambda d+4\lambda d^2)^{n+1}
\end{align*}
for \(x\neq 0\) and
\begin{align*}
|Q^{n+1}(0,y)|&=|QQ^n(0,y)|\\
&=|(1-4\lambda d)Q^n(0,y)+2\lambda\sum_{z:z\sim 0}Q^n(z,y)\\
&+\lambda\sum_{z:z\sim 0}\sum_{w:w\sim 0}Q(w+z,y)|\\
&\leq (1+4\lambda d+4\lambda d+4\lambda d^2)(1+8\lambda d+4\lambda d^2)^{n}\\
&=(1+8\lambda d+4\lambda d^2)^{n+1}.
\end{align*}
Therefore \eqref{poly bound of Q} holds for any \(n\geq 1\) by induction.

\qed

\begin{theorem}\label{unique solution of linear ode}
\(Q\) is the same as that in \eqref{Q}. There exists a unique solution \(\{f_t\in L^{\infty}(\textbf{Z}^d), t\geq0\}\) to the
following ODE
\[
\frac{d}{dt}f_t=Qf_t
\]
with initial condition that \(f_0\in L^{\infty}(\textbf{Z}^d)\). Moreover, \(f_t\) is with the form
\[
f_t=\exp\{tQ\}f_0
\]
and satisfies
\[
\|f_t\|_\infty\leq \exp\{1+8\lambda d+4\lambda d^2\}\|f_0\|_\infty
\] for any \(t\geq0\).
\end{theorem}
\proof
By a calculation similar with that in the proof of Theorem \ref{upper bound of Qelement}, it is easy to verify that
\[
\|Qf\|_\infty\leq(1+8\lambda d+4\lambda d^2)\|f\|_\infty
\]
for any \(f\in L^{\infty}(\textbf{Z}^d)\).
Therefore the linear operator \(T:f\rightarrow Qf\)  on \(L^{\infty}(\textbf{Z}^d)\) with norm \(\|\cdot\|_\infty\) satisfies the
Lipshitz condition. Then Theorem \ref{unique solution of linear ode} follows the classical theory of linear ODE on Banach spaces.

\qed

\begin{theorem}\label{positive of Qelement}
\(\exp\{tQ\}(x,y)\geq0\) for any \(x,y\in \textbf{Z}^d\).
\end{theorem}
\proof
We denote by \(I_{\textbf{Z}^d}\) the identity matrix \(\{\delta(x,y)\}_{x,y\in \textbf{Z}^d}\). Let
\[B=Q+4\lambda d I_{\textbf{Z}^d},\]
then \(B(x,y)\geq0\) for any \(x,y\in \textbf{Z}^d\) according to the definition of \(Q\). As a result, \(\exp\{tB\}(x,y)\geq0\) for any \(x,y\in \textbf{Z}^d\).
Since \(I_{\textbf{Z}^d}B=BI_{\textbf{Z}^d}\),
\[
\exp\{tQ\}=\exp\{tB\}\exp\{-4t\lambda dI_{\textbf{Z}^d}\}
\]
and hence
\[
\exp\{tQ\}(x,y)=\exp\{-4t\lambda d\}\exp\{tB\}(x,y)\geq0.
\]

\qed

\textbf{Acknowledgments.} This work is supported by the National Basic Research Program of China (2011CB808000), National Natural Science Foundation of China (No. 11001004) and China Scholarship Council (No. 201206010097).

{}
\end{document}